\documentclass{article}

\usepackage{amsmath}
\usepackage{amsthm}
\usepackage{amssymb}
\usepackage{amsfonts}
\usepackage{bm}
\usepackage{color}
\usepackage{comment}
\usepackage{enumerate}
\usepackage{float}
\usepackage{framed}
\usepackage{graphicx}
\usepackage{hyperref}
\usepackage[utf8]{inputenc}
\usepackage{lipsum} 
\usepackage{mathrsfs}
\usepackage{mathtools}
\usepackage{paralist}
\usepackage{url}
\usepackage{xargs}                      
\usepackage{xcolor}  
\usepackage[all,pdf]{xy}
\usepackage{stmaryrd}  

\newtheorem{theorem}{Theorem}
\newtheorem{lemma}[theorem]{Lemma}

\newtheorem{proposition}[theorem]{Proposition}
\newtheorem{definition}[theorem]{Definition}

\newcommand{\alg}{\mathrm{alg}}

\newcommand{\Id}{\mathrm{Id}}

\newcommand*\norm[1]{\left\lVert#1\right\rVert}
\newcommand*\inner[2]{\langle#1|#2\rangle}

\newcommand{\NN}{\mathbb{N}}
\newcommand{\ZZ}{\mathbb{Z}}

\newcommand{\RR}{\mathbb{R}}
\newcommand{\CC}{\mathbb{C}}

\begin{document}
	\title{A Positivstellensatz on the Matrix Algebra of Finitely Generated Free Group}
	
	\author{Hao Liang}
	
	
	
	\maketitle
	
	\begin{abstract}
		Positivstellens{\"a}tze are a group of theorems on the positivity of involution algebras over $\mathbb{R}$ or $\mathbb{C}$. One of the most well-known Positivstellensatz is the solution to Hilbert's 17th problem given by E. Artin, which asserts that a real polynomial in $n$ commutative variables is nonnegative on real affine space if and only if it is a sum of fractional squares. Let $m$ and $n$ be two positive integers. For the free group $F_n$ generated by $n$ letters, and a symmetric polynomial $b$ with variables in $F_n$ and with $n$-by-$n$ complex matrices coefficients, we use real algebraic geometry to give a new proof showing that $b$ is a sum of Hermitian squares if and only if $b$ is mapped to a positive semidefinite matrix under any finitely dimensional unitary representation of $F_n$. 
	\end{abstract}

	\section{Introduction}
	Positivstellens{\"a}tze are a group of theorems on the positive elements in an involution algebras (abbr., $*$-algebra) over $\RR$ and $\CC$. One of the most well-known Positivstellensatz is the solution to Hilbert's 17th problem given by E. Artin, which asserts that an $\RR$-polynomial in $n$ commutative variables is nonnegative on $\RR^n$ if and only if it is a sum of fractional squares \cite{Bochnak&Coste&Roy_1998}. For commutative algebras, several other Positivstellens{\"a}tze are found by Krivine \cite{Krivine1964}, Stengle \cite{Stengle1974}, Schm{\"u}dgen \cite{Schmuedgen1991} and Putinar \cite{Putinar1993}, see \cite{Marshall_2008} for an introduction. For the noncommutative case, the polynomial algebra and the points of $\RR^d$ are respectively substituted with a finitely generated $*$-algebra and a distinguished family of irreducible $*$-representations \cite{schmudgen2009noncommutative}. Various Positivstellens{\"a}tze are also summarized in \cite{schmudgen2009noncommutative,Netzer&Thom2013}.
	
	
	Let $m$ and $n$ be two positive integers. For the free group $F_n$ generated by $n$ letters, we denote by $\mathcal{A}:=\CC[F_n]$ the group algebra of $F_n$ over $\CC$. We consider the involution $*$ on the tensor algebra 
	$$M_m(\mathcal{A}):=M_m(\CC[F_n])=M_m(\CC)\otimes_{\CC}\CC[F_n]$$ 
	defined by $(M\otimes g)^*=M^{\dagger}\otimes g^{-1}$, where $M\in M_n(\CC)$ and $g\in F_n$. Using this involution, an unitary representation of $F_n$ is equivalent to a $*$-representation of $\CC[F_n]$ and $M_m(\CC[F_n])$. An element $b$ in $M_m(\mathcal{A})$ fixed by the involution $*$ is called a symmetric polynomial with variables in $F_n$ and with coefficients in $M_n(\CC)$. The set of symmetric elements in $M_m(\mathcal{A})$ is denoted by $M_m(\mathcal{A})^h$. A finite-dimensional $*$-representation of $\mathcal{A}$ is a $\CC$-algebraic homomorphism $\pi:\mathcal{A}\rightarrow M_N(\CC)$ preserving involution, where $N$ is a positive integer. In this article we apply the method in \cite[Theorem 6.1]{Netzer&Thom2013}\cite{Alekseev&Netzer&Thom_2019} to prove that $b$ is a sum of Hermitian squares, i.e., $b=\sum_{k=1}^{n}{a_k^* a_k}$ for some $a_k\in M_m(\mathcal{A})$, if and only if $b$ is mapped to a positive semidefinite matrix under any finitely dimensional unitary representation of $F_n$. 
	
	\begin{theorem}\label{MainTheorem}
		If $b=(b_{ij})\in M_m(\CC[F_n])^h$ is mapped to a positive semidefinite matrix $\pi_m(b):=(\pi(b_{ij}))\in M_{mN}(\CC)$ under every finite-dimensional unitary representation $\pi:F_n\rightarrow M_N(\CC)$, then $b\in\sum^h{M_n(\mathcal{A})}^2$.
	\end{theorem}
	
	Theorem \ref{MainTheorem} is proved by contradiction as follows. Using the separation theorem over a real closed field $\mathbf{R}$ proposed in \cite{Netzer&Thom2013}, there is a $\mathbf{C}$-state $f$ on $M_n(\mathcal{A}_{\mathbf{C}}):=M_n(\mathcal{A}\otimes_{\CC}\mathbf{C})$ strictly separating $b$ and the cone $\sum^h{M_n(\mathcal{A}_{\mathbf{C}})}^2$ of sums of Hermitian squares, where $\mathbf{C}=\mathbf{R}(\sqrt{-1})$ is the algebraic closure of $\mathbf{R}$. Then the GNS-construction for the $\mathbf{C}$-state $f$ gives a $*$-representation of $M_n(\mathcal{A}_{\mathbf{C}})$ on a preHilbert space over $\mathbf{C}$ such that $b$ is not sent to a positive operator. Using Choi's matrix trick \cite[Theorem 7]{Choi_1980} we can further construct a finite $*$-representation over $\mathbf{C}$ such that the image of $b$ is still not positive. Finally, Tarski's Transfer Principle transfers this finite representation from over $\mathbf{C}$ to over $\CC$.
	
	Our result generalizes the results in the previous literature \cite[Theorem 6.1]{Netzer&Thom2013} and \cite[Theorem 7.1]{Bakonyi&Timotin2007}. \cite[Theorem 6.1]{Netzer&Thom2013} proves that the case $m=1$ holds, and \cite[Theorem 7.1]{Bakonyi&Timotin2007} shows that $b$ is a sum of Hermitian squares if and only if $b$ is mapped to a positive semidefinite matrix under all unitary representation of $F_n$, and this result is also reproved in \cite{Ozawa2013}. 

	
	
	\section{Preliminaries}
	In this section we list necessary results for the proof of Theorem \ref{MainTheorem}.
	
	\subsection{Real closed field and Tarski's Transfer Principle}
	
	The reference for this section is \cite{Prestel&Delzell_2001}. An ordered field $(F,\leqslant)$ is a field $F$ equipped with a linear ordering ``$\leqslant$" satisfying
	\begin{enumerate}
		\item $\forall a,b,c\in F,(a\leqslant b\Rightarrow a+c\leqslant b+c)$, and
		
		\item $\forall a,b\in F,((a\geqslant0\wedge b\geqslant0)\Rightarrow ab\geqslant0)$.
	\end{enumerate}
	A real closed field $\mathbf{R}$ is a field such that $\mathbf{R}$ can be ordered and any proper algebraic field extension of $\mathbf{R}$ cannot be ordered, e.g., real algebraic number field $\RR_{\alg}$, real number field $\RR$ and real Puiseux series field $\RR(X)^{\wedge}$ (see \cite[Example 1.2.3]{Bochnak&Coste&Roy_1998}). Artin and Schreier have given a characterization for real closed field as follows 
	\begin{theorem}\cite{Artin&Schreier_1926}
		For a field $\mathbf{R}$ the following are equivalent
		\begin{enumerate}
			\item $\mathbf{R}$ is real closed;
			
			\item $\mathbf{R}$ can be (uniquely) ordered by $a\leqslant b:\Leftrightarrow\sqrt{b-a}\in\mathbf{R}$, and for any univariate polynomial $p\in\mathbf{R}[T]$ with odd degree, $p$ has a root in $\mathbf{R}$.
			
			\item $\mathbf{R}\subsetneq\mathbf{R}[\sqrt{-1}]$ and $\mathbf{R}[\sqrt{-1}]$ is algebraically closed;
		\end{enumerate}		
	\end{theorem}
	In this case we denote by $\mathbf{R}_{\geqslant0}$ the subset of elements not less than $0$ in $\mathbf{R}$ with respect to the unique ordering. A real closed field shares many similar algebraic-geometric properties with the real number field $\RR$, for example, Tarski's Transfer Principle, or equivalently Elimination of Quantifiers \cite{Prestel&Delzell_2001}.
	
	Given a subring $A$ of a real closed field $\mathbf{R}$ and a natural number $n\in\NN$, an $n$-ary semialgebraic definition $\delta(X_1,\dots,X_n)$ over $A$ is a Boolean combination (finite combinations of conjunction $\wedge$, disjunction $\vee$ and negation $\neg$) of formulae of the form $f(X_1,\dots,X_n)>0$, with $f\in A[X_1,\dots,X_n]$, e.g., for $A=\ZZ\subset\RR$ and $n=4$ let $\delta(X_1,\dots,X_4)$ be
	$$((X_2 X_1^2+X_3 X_1+X_4>0)\wedge(X_2>0))\vee(\neg(X_4<0)).$$
	For any real closed field $\mathbf{R}'$ containing $A$ such that $A\cap\mathbf{R}_{\geqslant0}=A\cap\mathbf{R}'_{\geqslant0}$, we can evaluate the variables $X_1,\dots,X_n$ at any point $p\in(\mathbf{R}')^n$ and check the truth value of $\delta(p)$, e.g., let $p=(1,1,1,1)$ then $\delta(p)=1$ is true in the above example, while $p=(1,-1,1,-1)$ implies $\delta(p)=0$ is false. Given such an $n$-ary semialgebraic definition $\delta=\delta(X_1,\dots,X_n)$ over $A$ we denote 
	$$\delta[\mathbf{R}']=\{p\in(\mathbf{R}')^n:\delta(p)~\text{is true}\}.$$
	
	Given a subring $A$ of a real closed field $\mathbf{R}$ and two natural number $m,n\in\NN$, an $n$-ary prenex definition $\sigma(X_1,\dots,X_n)$ over $A$ with $m$ constrained variables $Y_1,\dots,Y_m$ is like 
	$$Q_1 Y_1, Q_2 Y_2, \dots, Q_m Y_m, \delta(X_1,\dots,X_n,Y_1,\dots,Y_m),$$
	where $Q_i$ is $\exists$ or $\forall$, and $\delta$ is an $(m+n)$-ary semialgebraic definition over $A$. For any real closed field $\mathbf{R}'$ containing $A$ such that $A\cap\mathbf{R}_{\geqslant0}=A\cap\mathbf{R}'_{\geqslant0}$, we can also evaluate the variables $X_1,\dots,X_n$ at any point $p\in(\mathbf{R}')^n$ and check the truth value of $\sigma(p)$ over $\mathbf{R}'$, namely,
	$$\sigma(p)=\left[Q_1 y_1\in\mathbf{R}', Q_2 y_2\in\mathbf{R}', \dots, Q_m y_m\in\mathbf{R}', \delta(p_1,\dots,p_n,y_1,\dots,y_m)\right].$$
	We also denote
	$$\sigma[\mathbf{R}']=\{p\in(\mathbf{R}')^n:\sigma(p)~\text{is true}\}.$$
	
	\begin{theorem}\cite[Theorem 2.1.10, Tarski's Transfer Principle]{Prestel&Delzell_2001}
		
		Let $\mathbf{R}_1,\mathbf{R}_2$ be two real closed fields inducing the same ordering on a common subring $A$, and $\sigma(\underline{X}_n)=Q_1 Y_1,\dots,Q_m Y_m,\delta(\underline{X}_n,\underline{Y}_m)$ is an $n$-ary prenex definition over $A$ for $m,n\in\NN$. Given $c=(c_1,\dots,c_n)\in A^n$, then
		\begin{align*}
			&Q_1 y_1\in\mathbf{R}_1,\dots,Q_m y_m\in\mathbf{R}_1,\delta(c,\underline{y}_m),~\text{and}\\
			&Q_1 y_1\in\mathbf{R}_2,\dots,Q_m y_m\in\mathbf{R}_2,\delta(c,\underline{y}_m) 
		\end{align*}
		has the same truth value. 
	\end{theorem}
	
	\subsection{$*$-rings, quadratic modules and $*$-representations}
	
	A $*$-ring is an associative and unital ring $\mathcal{A}$ equipped with an involution $*$, i.e., for $a,b,c\in\mathcal{A}$ we have 
	$$(a+b)^*=a^*+b^*,~(ab)^*=b^*a^*,~(a^*)^*=a.$$ 
	Note that $0^*=0$ and $1^*=1$, so the restriction of $*$ on $\ZZ\cdot1$ is trivial. We also remark that the center $Z(\mathcal{A})$ of a $*$-ring is again a $*$-ring. The set $\mathcal{A}^h$ of $*$-invariant element of $\mathcal{A}$ is an Abelian subgroup of the additive group $\mathcal{A}$, called the symmetric subgroup of $\mathcal{A}$. A $*$-homomorphism between two $*$-ring $(\mathcal{A},*)$ and $(\mathcal{B},*)$ is a ring map $\varphi:\mathcal{A}\rightarrow\mathcal{B}$ preserving $*$; i.e., $\varphi(a^*)=\varphi(a)^*$. 
	
	
	Let $\mathbf{R}$ be a real closed field and $\mathbf{C}=\mathbf{R}(\sqrt{-1})$ be the algebraic closure of $\mathbf{R}$. A complex closed field is the field extension $\mathbf{R}\subseteq\mathbf{C}$, or equivalently the $*$-ring $(\mathbf{C},*)$ with the conjugate involution $*$ over $\mathbf{R}$. A $(\mathbf{C},*)$-algebra is a $*$-homomorphism $\varphi:(\mathbf{C},*)\rightarrow(\mathcal{A},*)$ such that $\varphi(\mathbf{C})$ is contained in the center $Z(\mathcal{A})$ of $\mathcal{A}$. The symmetric subgroup $\mathcal{A}^h$ is an $\mathbf{R}$-linear subspace of $\mathcal{A}$, and we have $\mathcal{A}=\mathcal{A}^h\oplus\sqrt{-1}\mathcal{A}^h$. A $(\mathbf{C},*)$-map between two $(\mathbf{C},*)$-algebras is a $*$-homomorphism $\psi:(\mathcal{A},*)\rightarrow(\mathcal{B},*)$ which is also a $\mathbf{C}$-algebra map. Note that $\psi(\mathcal{A}^h)\subseteq\mathcal{B}^h$ and $\psi(\sqrt{-1}\mathcal{A}^h)\subseteq\sqrt{-1}\mathcal{B}^h$. 
	
	For a $*$-ring $A$, a subset $Q\subseteq\mathcal{A}^h$ is called a quadratic module if $Q+Q\subseteq Q,~1\in Q,~-1\notin Q$ and $a^* Q a\subseteq Q$ for any $a\in\mathcal{A}$. Any quadratic module contains $\sum^h{\mathcal{A}^2}$, the set of sums of Hermitian squares, and $\sum^h{\mathcal{A}^2}$ is the smallest quadratic module of $\mathcal{A}$. A quadratic module $Q$ of a $(\mathbf{C},*)$-algebra $(\mathcal{A},*)$ is called $(\mathbf{C},*)$-Archimedean if for any $a\in\mathcal{A}$ there is an $r\in \mathbf{R}_{\geqslant0}$ such that $r-a^* a\in Q$. A $(\mathbf{C},*)$-Archimedean quadratic module $Q$ of a $(\mathbf{C},*)$-algebra $(\mathcal{A},*)$ is called $(\mathbf{C},*)$-Archimedean closed if any $a\in\mathcal{A}^h$ satisfying $a+\varepsilon\in Q$ for any $\varepsilon\in\mathbf{R}_{>0}$ must lie in $Q$. 

	An inner product on a $\mathbf{C}$-linear space $\mathcal{H}$ is a positive sesquilinear form $$\inner{-}{-}:\mathcal{H}\times\mathcal{H}\rightarrow\mathbf{C},$$
	i.e., for $u,v,w\in\mathcal{H},~a,b\in\mathbf{C}$, we have $\inner{u+v}{w}=\inner{u}{w}+\inner{v}{w},~\inner{u}{v+w}=\inner{u}{v}+\inner{u}{w},~\inner{au}{bv}=\overline{a}b\inner{u}{v},~\inner{a}{b}=\overline{\inner{b}{a}}$, and $a\neq0\Rightarrow\inner{a}{a}\in\mathbf{R}_{>0}$. A $\mathbf{C}$-linear space $\mathcal{H}$ equipped with an $\mathbf{C}$-inner product $\inner{-}{-}$ is called a $(\mathbf{C},*)$-preHilbert space. A $\mathbf{C}$-linear map $T$ from $\mathcal{H}_1$ to $\mathcal{H}_2$ is called unitary if $T$ is bijective and isometric.
	
	For a complex closed field $(\mathbf{C},*)$, let $\mathcal{H}$ be a $(\mathbf{C},*)$-preHilbert space. Given a $\mathbf{C}$-linear operator $A$ on $\mathcal{H}$, the adjoint operator $A^*$ of $A$ is the (unique) operator on $\mathcal{H}$ satisfying 
	$$\inner{A^* u}{v}=\inner{u}{Av},~\forall u,v\in\mathcal{H}.$$
	Note that the adjoint operator $A^*$ of a given operator $A$ does not necessarily exist. However, if $\mathcal{H}$ is of finite dimension over $\mathbf{C}$ then every $\mathbf{C}$-linear operator on $\mathcal{H}$ must have an (unique) adjoint. We can straightforward check that
	\begin{enumerate}
		\item If $A^*$ exists, then so does $(A^*)^*$, and $(A^*)^*=A$, 
		
		\item If both $A^*$ and $B^*$ exist, then $(A+B)^*$ and $(AB)^*$ also exist, and $(A+B)^*=A^*+B^*$ as well as $(AB)^*=B^*A^*$, and
		
		\item $\lambda^*=\overline{\lambda}$ for $\lambda\in\mathbf{C}$.
	\end{enumerate}
	Therefore, the set of operators on $\mathcal{H}$ whose adjoint operator exists is a $(\mathbf{C},*)$-algebra, which is denoted by $(\mathcal{L}^*(\mathcal{H}),*)$, or simply $\mathcal{L}^*(\mathcal{H})$. Note that if $\mathcal{H}$ is of finite dimension over $\mathbf{C}$ then $\mathcal{L}^*(\mathcal{H})$ coincides with the $\mathbf{C}$-linear operator algebra on $\mathcal{H}$. A $\mathbf{C}$-linear operator $A$ on $\mathcal{H}$ is called positive if $A\in\mathcal{L}^*(\mathcal{H})^h$ and $\inner{u}{Au}\geqslant0$ for all $u\in\mathcal{H}$. 
	
	A $(\mathbf{C},*)$-algebra map $(\mathcal{A},*)\rightarrow(\mathcal{L}^*(\mathcal{H}),*)$ for some $(\mathbf{C},*)$-preHilbert space $\mathcal{H}$ is also called a $(\mathbf{C},*)$-representation of $(\mathcal{A},*)$, or simply a $*$-representation of $(\mathcal{A},*)$. If moreover $\mathcal{H}$ is of finite dimension over $\mathbf{C}$, then this $(\mathbf{C},*)$-representation is called finite. Given two $(\mathbf{C},*)$-representations $\pi_1$ and $\pi_2$ of $(\mathcal{A},*)$ on $(\mathbf{C},*)$-preHilbert space $\mathcal{H}_1$ and $\mathcal{H}_2$ respectively, if there exists a unitary map $T:\mathcal{H}_1\rightarrow\mathcal{H}_2$ such that for all $a\in\mathcal{A}$ we have $T\circ\pi_1(a)=\pi_2(a)\circ T$, then $\pi_1$ and $\pi_2$ are called unitarily equivalent. 
	
	For a complex closed field $(\mathbf{C},*)$ and a group $G$, we write $\mathbf{C}[G]$ for the group algebra of the group $G$ over $\mathbf{C}$, and we equip $\mathbf{C}[G]$ with the involution determined by $g^*=g^{-1}$ and $\lambda^*=\overline{\lambda}$ for $g\in G$ and $\lambda\in\mathbf{C}$. In this way, $(\mathbf{C}[G],*)$ becomes a $(\mathbf{C},*)$-algebra. We denote the $*$-invariant subspace of $\mathbf{C}[G]$ by $\mathbf{C}[G]^h$. A unitary representaion of $G$ over the complex closed field $(\mathbf{C},*)$ is a $(\mathbf{C},*)$-algebra map from $(\mathbf{C}[G],*)$ to $(\mathcal{L}^*(\mathcal{H}),*)$ for some $(\mathbf{C},*)$-preHilbert space. 
	
	A $*$-ring $(\mathcal{A},*)$ is called real reduced if for any elements $a_1,\dots,a_n\in\mathcal{A}$, $\sum_{i}{a_i^* a_i}=0$ implies all the $a_i$'s are $0$ (see \cite{Netzer&Thom2013}). For the group $(\mathbf{C},*)$-algebra $\mathbf{C}[G]$ of the group $G$ over the complex closed field $\mathbf{C}$, using the faithful trace
	\begin{align*}
		\tau:\mathbf{C}[G]&\rightarrow\mathbf{C}\\
		\sum_{g\in G}{a_g g}&\mapsto a_1
	\end{align*}
	\cite{Netzer&Thom2013} assserts that $\mathbf{C}[G]$ is real reduced. 
	
	\subsection{States and GNS-construction of $(\mathbf{C},*)$-algebras}
	
	Still let $\mathbf{R}$ be a real closed field and $\mathbf{C}=\mathbf{R}(\sqrt{-1})$ be the algebraic closure of $\mathbf{R}$. A $(\mathbf{C},*)$-state of a $(\mathbf{C},*)$-algebra $(\mathcal{A},*)$ is a $\mathbf{C}$-linear function $f:\mathcal{A}\rightarrow\mathbf{C}$ such that $f(1)=1$ and $f(a^* a)\in \mathbf{R}_{\geqslant0}$ for all $a\in\mathcal{A}$. Given a $(\mathbf{C},*)$-representation $\pi:(\mathcal{A},*)\rightarrow(\mathcal{L}^*(\mathcal{H}),*)$ of $(\mathbf{C},*)$-algebra $(\mathcal{A},*)$, for any unit vector $\xi\in\mathcal{H}$, we have an associated state $f_{\pi,\xi}:\mathcal{A}\rightarrow \mathbf{C}$ sending $a\in\mathcal{A}$ to $\inner{\xi}{\pi(a)\xi}$. Conversely, 
	
	\begin{theorem}[Gelfand-Nairmark-Segal construction]
		
		Given a $(\mathbf{C},*)$-state $f$ of a $(\mathbf{C},*)$-algebra $(\mathcal{A},*)$, there is a unique cyclic $(\mathbf{C},*)$-representation $\pi_f$ (up to isomorphism) of $(\mathcal{A},*)$ on a $(\mathbf{C},*)$-preHilbert space $\mathcal{H}_f$ with a unit vector $\xi_f\in\mathcal{H}_f$, namely, $\mathcal{H}_f=\mathcal{A}/N_f$ for
		$$N_f=\{a\in\mathcal{A}:f(a^* a)=0\},$$
		$\inner{a+N_f}{b+N_f}_{\mathcal{H}_f}=f(a^* b)$, $\pi_f(a)(b+N_f)=ab+N_f$ and $\pi_f(a)^*=\pi_f(a^*)$ for all $a,b\in\mathcal{A}$ and $\xi_f=1+N_f$, such that $f(a)=\inner{\xi_f}{\pi_f(a)\xi_f}_{\mathcal{H}_f}$ for $a\in\mathcal{A}$. 
	\end{theorem}
	
	A $(\mathbf{C},*)$-state $f$ of a $(\mathbf{C},*)$-algebra $(\mathcal{A},*)$ is called finite if the associated $(\mathbf{C},*)$-preHilbert space $\mathcal{H}_f$ is of finite dimension over $\mathbf{C}$. 
	
	For the Archimedean case, the following proposition provides a positive functional strictly separating sums of Hermitian squares and any other point.
	
	\begin{proposition}
		A $(\CC,*)$-Archimedean quadratic module $Q$ of a $(\CC,*)$-algebra $(\mathcal{A},*)$ is $(\CC,*)$-Archimedean closed if and only if for any $a\in\mathcal{A}^h\setminus Q$ there is a $(\CC,*)$-state $f$ of $(\mathcal{A},*)$ such that $f(a)<0$ and $f(Q)\geqslant0$. 
	\end{proposition}
	
	\begin{proof}
		``if'': Suppose that for any $a\in\mathcal{A}^h$ not lying in $Q$ there is a $(\CC,*)$-state $f$ of $(\mathcal{A},*)$ such that $f(a)<0$ and $f(Q)\geqslant0$. Then we have $$Q=\bigcap\{f^{-1}(\RR_{\geqslant0}):\text{$f$ is a state, and $f(Q)\geqslant0$}\},$$
		which implies that $Q$ is $(\CC,*)$-Archimedean closed.
		
		``only if": Assume that $Q$ is $(\CC,*)$-Archimedean closed. Then by Eidelheit-Kakutani Separation Theorem \cite{Barvinok_2002}, for any $a\in\mathcal{A}^h$ not lying in $Q$ there is a $\RR$-linear function $f:\mathcal{A}^h\rightarrow\RR$ such that $f(a)<0$ and $f(Q)\geqslant0$. If $f(1)=0$ then by Cauchy-Schwartz inequality we will obtain a contradiction that $f=0$. Hence $f(1)>0$, and we may rescale $f$ and assume $f(1)=1,~f(Q)\geqslant0$ and $f(a)<0$. In this way we find a desired state $f$.
	\end{proof}
	
	For a real reduced $(\CC,*)$-algebra, \cite[Lemma 3.6 and Theorem 3.11]{Netzer&Thom2013} provides a weaker separation theorem. 
	
	\begin{proposition}\cite[Lemma 3.6 and Theorem 3.11]{Netzer&Thom2013}
		
		Let $(\mathcal{A},*)$ be a real reduced $(\CC,*)$-algebra, and $b\in\mathcal{A}^h\setminus\sum^h{\mathcal{A}^2}$. Then there is real closed extension field $\mathbf{R}$ of $\RR$ and a $(\mathbf{C},*)$-state $f:\mathcal{A}\otimes_{\CC}\mathbf{C}\rightarrow\mathbf{C}$ for $\mathbf{C}=\mathbf{R}(\sqrt{-1})$ such that $f(b)< 0$ and $f(a^*a)>0$ for all $a\in A\setminus0$.
	\end{proposition}
	
	\subsection{Tensor of $(\mathbf{C},*)$-algebras}
	
	\begin{definition}
		Let $(\mathbf{C},*)$ be a complex closed field. The tensor $(\mathbf{C},*)$-algebras of two $(\mathbf{C},*)$-algebras $(\mathcal{A},*)$ and $(\mathcal{B},*)$ is the $(\mathbf{C},*)$-algebra $(\mathcal{A}\otimes_{\mathbf{C}}\mathcal{B},*)$ where the involution is given by
		$$\left(\sum_k{a_k\otimes b_k}\right)^*=\sum_k{a_k^*\otimes b_k^*}.$$
	\end{definition}
	
	Recall that, for any positive integer $n$ we have a matrix $(\mathbf{C},*)$-algebra $(M_n(\mathbf{C}),*)$, where $A^*=A^{\dagger}$, i.e., $(a_{ij})^*=(\overline{a_{ji}})$. Provided a $(\mathbf{C},*)$-algebra $(\mathcal{A},*)$, using the $(\mathbf{C},*)$-isomorphism
	\begin{align*}
		M_n(\mathcal{A})&\rightarrow M_n(\mathbf{C})\otimes_{\mathbf{C}}\mathcal{A},\\
		(a_{ij})&\mapsto\sum_{i,j}{E_{ij}\otimes a_{ij}}
	\end{align*}
	$M_n(\mathcal{A})$ becomes a $(\mathbf{C},*)$-algebra, and the involution is the natural one
	$$(a_{ij})^*:=(a_{ji}^*).$$
	
	Given a $(\mathbf{C},*)$-representation $\pi:\mathcal{A}\rightarrow\mathcal{L}^*(\mathcal{H})$ of the $(\mathbf{C},*)$-algebra $\mathcal{A}$, there is a natural $(\mathbf{C},*)$-representation $\pi_n$ of the $(\mathbf{C},*)$-algebra $M_n(\mathcal{A})$ on $\mathcal{H}^n$ satisfying 
	$$\pi_n((a_{ij}))=(\pi(a_{ij})).$$
	Conversely, provided a $(\mathbf{C},*)$-representation $\sigma:M_n(\mathcal{A})\rightarrow\mathcal{L}^*(\mathcal{K})$, there exist a $(\mathbf{C},*)$-preHilbert space $\mathcal{H}$ and a $(\mathbf{C},*)$-representation $\pi:\mathcal{A}\rightarrow\mathcal{L}^*(\mathcal{H})$ such that $\mathcal{K}=\mathcal{H}^n$ and $\sigma((a_{ij}))=(\pi(a_{ij}))$.
	
	Recall that a $*$-ring $(\mathcal{A},*)$ is called real reduced if for any finite elements $a_1,\dots,a_n\in\mathcal{A}$, $\sum_{i}{a_i^* a_i}=0$ implies all the $a_i$'s are zero \cite{Netzer&Thom2013}. 
	
	\begin{lemma}\label{lemma7}
		A $(\mathbf{C},*)$-algebra $(\mathcal{A},*)$ is real reduced if and only if $(M_n(\mathcal{A}),*)$ is real reduced.
	\end{lemma}
	
	\begin{proof}
		``If'': Suppose that $(M_n(\mathcal{A}),*)$ is real reduced. If $a_1^* a_1+\dots+a_m^* a_m=0$ in $\mathcal{A}$, then by embedding $\mathcal{A}$ into the $(1,1)$-entry of $M_n(\mathcal{A})$ we see that $a_i$'s are all zero. Hence $(\mathcal{A},*)$ is real reduced.
		
		``Only if'': Suppose that $(\mathcal{A},*)$ is real reduced, and $A:=A_1^* A_1+\dots+A_m^* A_m=0$ for $A_k=(a_{ij,k})$ in $M_n(\mathcal{A})$. Notice that the $(l,l)$-entry of $A$ is $\sum_{i,k}{a_{il,k}^* a_{il,k}}$, hence from $A=0$ we see that all $A_k$'s are zero.
	\end{proof}

	\section{Proof of Theorem \ref{MainTheorem}}
		We can now give the proof for Theorem \ref{MainTheorem}.
		
		\begin{proof}[Proof of Theorem \ref{MainTheorem}]
			Assume $b=(b_{ij})\notin\sum{M_n(\mathcal{A})}^2$. Since $\mathcal{A}$ is real reduced \cite[Remark 4.4.]{Netzer&Thom2013}, Lemma \ref{lemma7} implies that $M_n(\mathcal{A})$ is also real reduced. According to \cite{Netzer&Thom2013}, there is a complex closed field extension $(\mathbf{C},*)$ over $(\CC,*)$ and a $(\mathbf{C},*)$-state $f:\mathbf{C}\otimes_{\CC}M_n(\mathcal{A})\rightarrow\mathbf{C}$ such that $f(b)<0$. The GNS construction associated to $f$ gives a cyclic $(\mathbf{C},*)$-representation $(\pi_f,\xi_f)$ of $\mathbf{C}\otimes_{\CC}M_n(\mathcal{A})=M_n(\mathbf{C}\otimes_{\CC}\mathcal{A})$ such that $\inner{\xi_f}{\pi_f(b)\xi_f}<0$. Write $\pi_f=M_n(\mathbf{C})\otimes_{\CC}\pi$ for $\pi:\mathbf{C}\otimes_{\CC}\mathcal{A}\rightarrow\mathcal{L}^*(\mathcal{H})$, then $\mathcal{H}_f=\mathcal{H}^n$. Denote $\xi_i=E_{ii}\xi_f\in\mathcal{H}$ for $i\in[n]$, then $\xi_f=(\xi_1,\dots,\xi_n)\in\mathcal{H}^n,~\norm{\xi_f}^2=\sum_i{\norm{\xi_i}^2}=1$ and $\inner{\xi_f}{\pi_f(b)\xi_f}=\sum_{i,j}{\inner{\xi_i}{\pi(b_{ij})\xi_j}}$. 
			
			Now let $\mathcal{H}'$ be the finite-dimensional $(\mathbf{C},*)$-subspace of $\mathcal{H}$, generated by $\pi(w)\xi_k$ for $k\in[n]$ and all words $w\in F_n$ of length at most $d$, where $d$ is the maximal word length in $b$. Using the usual Gram–Schmidt procedure over $\mathbf{C}$, we can find an orthogonal projection map $p:\mathcal{H}\rightarrow\mathcal{H}'$ such that for the natural inclusion $q:\mathcal{H}'\rightarrow\mathcal{H}$ we have $pq=\Id_{\mathcal{H}'}$ and $(qp)^2=qp\in\mathcal{L}^*(\mathcal{H})$. Define 
			$$C_i=p\circ\pi(X_i)\circ q\in\mathcal{L}(\mathcal{H}')=\mathcal{L}^*(\mathcal{H}').$$
			Since all $C_i$'s are contractions, we see that the linear operators $\sqrt{1-C_i^* C_i}$ and $\sqrt{1-C_i C_i^*}$ exist on $\mathcal{H}'$. Moreover, using polar decomposition we have 
			$$C_i\sqrt{1-C_i^* C_i}=\sqrt{1-C_i C_i^*}C_i.$$ 
			Hence the operators
			\begin{equation*}
				U_i=
				\begin{pmatrix}
					C_i & -\sqrt{1-C_i C_i^*}\\
					\sqrt{1-C_i^* C_i} & C_i^*
				\end{pmatrix}
				\in\mathcal{L}(\mathcal{H}'\oplus\mathcal{H}')
			\end{equation*}
			are unitary, and we can define a $(\mathbf{C},*)$-representation $\widetilde{\pi}$ of $C[F_n]$ on $\mathcal{H}'\oplus\mathcal{H}'$ by sending $X_i$ to $U_i$. Write $\xi_i'$ for the vector in $(\mathcal{H}'\oplus\mathcal{H}')^n$ whose $i$-th entry is $(\xi_i,0)$ and all other entries are zero, then we can check that $\inner{\xi_i'}{\widetilde{\pi}(b_{ij})\xi_j'}=\inner{\xi_i}{\pi(b_{ij})\xi_j}$. 
			
			Let $\widetilde{\pi}_n$ be the $(\mathbf{C},*)$-representation of $M_n(\mathbf{C}[F_n])$ on $(\mathcal{H}'\oplus\mathcal{H}')^n$ that maps $(a_{ij})\in M_n(\mathbf{C}[F_n])$ to $(\widetilde{\pi}(a_{ij}))\in\mathcal{L}((\mathcal{H}'\oplus\mathcal{H}')^n)$. For the vector $\xi'=(\xi_1',\dots,\xi_n')$ we see that $\inner{\xi'}{\widetilde{\pi}_n(b)\xi'}=\inner{\xi}{\pi_f(b)\xi}<0$. Using Tarski's Transfer Principle, we deduce that there is a finite-dimensional $(\CC,*)$-representation $\rho$ of $M_n(\CC[F_n])$ such that $\rho(b)$ is not positive semidefinite, which contradicts to the condition that $b$ is positive. Therefore, $b\in\sum{M_n(\mathcal{A})}^2$. 
		\end{proof}

	\bibliographystyle{alpha}
	\bibliography{LiangBook,LiangArticle}
\end{document}